\documentclass[a4paper,11pt]{article}

\usepackage{geometry, amsmath,amssymb,latexsym,enumerate}
\usepackage{xypic}

\begin{document}

\title{ Double groupoids, matched pairs and then matched triples}
 \author{Ronald Brown  \\ School of Computer Science\\
Bangor University  \\Gwynedd LL57 1UT,  U.K.} \maketitle

\section{Introduction}
\label{sec-int} In this note we show that the known relation between
double groupoids and matched pairs of groups may be extended, or
seems to extend,  to the triple case. The references give some other
occurrences of double groupoids.

I hope that someone can pursue these ideas further.
\section{Double groupoids} A {\em double groupoid} is simply a
groupoid object internal to the category of groupoids. Thus it
consists of a set  $G$ with two groupoid  structures, say
$\circ_1, \circ_2$ which satisfy the {\em interchange law},
namely the condition that there is only one way of evaluating
the composition  \begin{equation} \left[ \begin{array}{cc} x & y
\\ z & w \end{array}  \right]  { \spreaddiagramrows{-1.2pc}
\spreaddiagramcolumns{-1.2pc} \diagram . \rto^>{2} \dto^>{1}& \\
&  \enddiagram} \end{equation}  This gives the formula on a line
\begin{equation} (x\circ_1z)\circ_2(y\circ_1w)  =
(x\circ_2y)\circ_1(z\circ_2 w), \end{equation}  for all $x,y,z,w
\in G$ such that all the compositions on both sides are defined.

\section{Matched pairs of groups}

The aim of this section is to give an exposition of matched pairs of
groups in a style in keeping with the geometry of double groupoids,
and where the notation is suggestive for higher dimensions. For
information on matched pairs, see \cite{Maj}. I think this
relationship is well known but currently do not have a reference.

Let the group $G$  be given as the product $MN$ of subgroups
$M,N$ such that $M \cap N ={1}.$ Then each element $g$ of $G$
can be uniquely written as $mn, m\in M, n\in N.$  By taking
inverses we see that $g$ can also be written uniquely as $n'm',
n'\in N, m' \in M.$ We therefore write for $m\in M, n \in N$
$$mn= {}^mnm^n.$$ We write this pictorially as $$\diagram .
\dto_m \rto^{{}^mn} & . \dto^{m^n} \\ . \rto_n &. \enddiagram $$
 The advantage of this approach is that this easily gives rules
for products and these operations, as follows.  \par Consider
the diagram: $$\diagram . \dto_m \rto^{{}^mn} & .
\dto^{m^n}\rto^{{}^{m^n}p}& . \dto^{m^{np}}  \\ . \rto_n &.
\rto_p &. \enddiagram  $$ We obtain immediately that: $$m^{np}=
(m^n)^p, \qquad {}^m(np) = {}^mn( {}^{m^n}p).$$  Putting
$p=n^{-1}$ in this equation gives $$({}^mn)^{-1} =
^{m^n}(n^{-1}).$$  Consider  the diagram: $$\diagram .\dto_l
\rto^{{}^{lm}n} &.  \dto^{l^{{}{^m}n}} \\ . \dto_m \rto^{{}^mn}
& . \dto^{m^n} \\ . \rto^n & . \enddiagram  $$ This gives the
rules: $${}^{lm}n = {}^l({}^mn), \qquad (lm)^n = (l^{{}^mn})
m^n .$$ In this case we deduce $$(m^n)^{-1} = (m^{-1})^{{}^mn}.
$$ We also note that $$n^{-1}m^{-1}=(m^n)^{-1} (^mn)^{-1}.$$ If
we replace $n^{-1}$ with $n$, and $m^{-1}$ with $m$, and write
$\bar{m}=m^{-1}, \bar{n} =n^{-1}$, then we deduce that
\begin{equation} nm = (\bar{m}^{\bar{n}})^{-1}
({}^{\bar{m}}\bar{n})^{-1}= (m^{{}^{\bar{m}}\bar{n}})
({}^{\bar{m}^{\bar{n}}}n). \label{invert} \end{equation}  We can
express this in the language of double groupoids by saying that
the set $M\times N$ may be given the two action groupoid
structures $M\rtimes N$ and $M\ltimes N$. That is, in accordance
with the pictures above, we have $$\begin{array}{ccccccc}  (m,n)
& \circ_2 & (m^n,p) & = & (m,np) & \in  & M\rtimes N, \\
(l,{}^mn) & \circ_1 & (m,n) & = & (lm,n) & \in & M \ltimes N .
\end{array} $$ This double groupoid has the special property
that given any two adjacent edges then there is a unique square
filling them. Such double groupoids are well known to give rise
to a groupoid, by filling in according to the following diagram:
$$
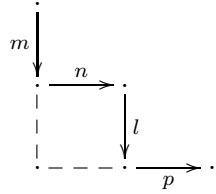
\diagram  . \dto_m & & \\ . \rto^n \xdashed[1,0] &. \dto^{l} &
\\ . \rdashed  & . \rto_{p} & . \enddiagram $$ leading to the
definition in our case that $$(m,n)(l,p) = (m
l^{{}^{\bar{l}}\bar{n}}, {}^{\bar{l}^{\bar{n}}}np).$$ \par

\section{Triple groupoids}  Here is an initial experiment on the
ideas. We have three groups $M,N,P$ for which \begin{itemize}
\item $M$ operates on the left of $N$ and of $P$ \item $N$
operates  on the right of $M$ and the left of $P$ \item $P$
operates on the right of $M$  and of $N$. \end{itemize} The
cubical model is then of the following form:
 $$
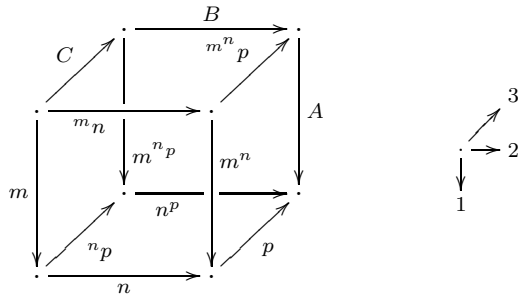
\diagram &  .\rrto^{B}\xline '[1,0] [2,0]^{m^{^np}}|>\tip &&.
\ddto^{A}\\. \urto^{C} \rrto_(0.3){^mn}\ddto_{m}&&.
\urto^{{}^{m^n}p} \ddto^(0.3){m^n} & \\&. \xline'[0,1]_{n^p}
[0,2]|>\tip && . \\ . \rrto_{n}\urto_{^np} &&. \urto_p &
\enddiagram \qquad   \qquad \spreaddiagramcolumns{-1.2pc}
\spreaddiagramrows{-1.2pc}\def\objectstyle{\ssize} \diagram &3
\\ . \dto  \rto \urto & 2\\ 1& \enddiagram$$ where \[
\begin{array}{cclcr} A & = &  m^{np}& = &(m^{^np})^{n^p}, \\ B&=
&  {}^{m^{^np}}(n^p)& = & (^mn)^{^{m^n}p},\\ C & =&  {}^{mn}p&
= &   {}^{^mn}(^{m^n}p). \end{array}   \]

At this stage one can begin to see the notation becoming
impossible to handle, for example to describe the reconstitution
of a group from the above data. We therefore adopt a \TeX\
style approach to superscripts and to left and right actions and
write $$^{\wedge}  mn = {}^mn, \qquad m^{\wedge}  n = m^n.$$ That is
a $^{\wedge}$ before a term raises it to a superscript. The above
picture then becomes as follows: $$\spreaddiagramrows{1.2pc}
\spreaddiagramcolumns{1.2pc}
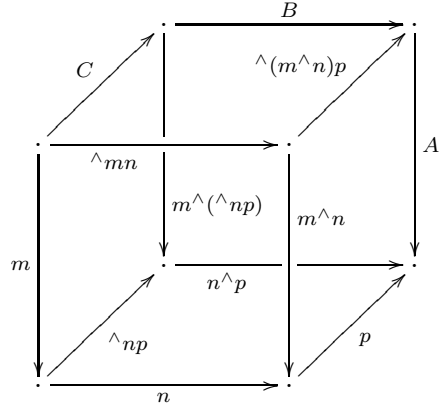
\diagram &  .\rrto^{B}\xline '[1,0]
[2,0]^{ m^{\wedge}  (^{\wedge}  np)}|>\tip &&. \ddto^{A}\\.
\urto^{C} \rrto_(0.3){^{\wedge}  mn}\ddto_{m}&&. \urto^{^{\wedge}
(m^{\wedge}  n)p} \ddto^(0.3){m^{\wedge}  n} & \\&.
\xline'[0,1]_{n^{\wedge}  p} [0,2]|>\tip && . \\ .
\rrto_{n}\urto_{^{\wedge}  np} &&. \urto_p & \enddiagram $$ where
\[ \begin{array}{cclcr} A & = &  m^{\wedge}  (np)& = &(m^{\wedge}
((^{\wedge}  n)p))^{\wedge}  (n^{\wedge}  p), \\ B&= &  ^{\wedge}
(m^{\wedge}  (^{\wedge} n p))(n^{\wedge}  p)& = & (^{\wedge}  mn)^{\wedge}
 (^{\wedge}  (m^{\wedge}  n)p),\\ C & =&  ^{\wedge}  (mn)p&  = &
^{\wedge}  (^{\wedge}  mn)(^{\wedge}  (m^{\wedge}  n)p). \end{array}
\]

This allows us to describe a possible formula for the
composition in the group associated with a matched triple.
First, the previous formula (\ref{invert}) for $nm$ may now be
written as \begin{equation}  nm  =  (m ^{\wedge} (^{\wedge}
\bar{m}\bar{n}))(^{\wedge} (^{\wedge} \bar{m} \bar{n})n).
\end{equation} Even this is too cumbersome, and we therefore
define two pairings by \begin{eqnarray} m\nearrow n & = & m
^{\wedge} (^{\wedge} \bar{m}\bar{n}), \\  m\nwarrow n & = & ^{\wedge}
(^{\wedge} \bar{m} \bar{n})n.  \end{eqnarray}

Now adjoining two cubes suggests the formula for a group
composition:$$(m,n,p)(\mu , \nu , \pi ) = (m((\mu \nearrow p)
\nearrow n), ((\mu \nearrow p) \nwarrow n)(\nu \nearrow (\mu
\nwarrow p)), (\nu \nearrow (\mu \nwarrow p)) \pi ). $$ There
remains to check that this satisfies the group axioms. This
should be done as far as possible not by manipulating these
formulae on a line but by examining the required cubical
diagrams.

This situation should occur for a matched triple of subgroups
$M,N,P$ of $G$ in which each subpair is a matched pair.

\end{document}